\documentclass{amsart}
\usepackage{graphicx}
\usepackage{fullpage,amsfonts,amsmath,amscd,amssymb}

\theoremstyle{plain}

\newtheorem*{lem}{Lemma}
\newtheorem*{prop}{Proposition}
\newtheorem*{thm}{Theorem}

\newtheorem*{cor}{Corollary}

\theoremstyle{remark}

\newtheorem*{rem}{Remark}
\newtheorem*{rems}{Remarks}

\newcommand{\abs}[1]{\left\vert#1\right\vert}

\newcommand{\C}{\mathbb C}
\newcommand{\Q}{\mathbb Q}

\newcommand{\To}{\longrightarrow}

\newcommand{\m}{\mathfrak{m}}
\newcommand{\p}{\mathfrak{p}}
\newcommand{\g}{\mathfrak{g}}
\newcommand{\eq}{\quad = \quad}
\newcommand{\gk}{\mathrm{GK}\!-\!\mathrm{dim}_k}
\newcommand{\Z}{\mathcal{Z}}
\begin{document}

\title[]{Unruffled extensions and flatness over central subalgebras}%
\author{Kenneth A. Brown}%
\address{Department of Mathematics, University of Glasgow, Glasgow G12 8QW}%
\email{kab@maths.gla.ac.uk}%

\thanks{Some of the research described here was
carried out during a workshop in Seattle in August 2003 funded by
Leverhulme Research Interchange Grant F/00158/X. Most of the paper
was written during a two-month stay at the Mittag-Leffler
Institute, Stockholm, in the autumn of 2003.
I am very grateful to the staff of the Institute for their hospitality. }%
\thanks{Extensive conversations with Ken Goodearl, Iain Gordon and
Toby Stafford have greatly helped me in this work>}
%\subjclass{}%

\date{March 22, 2004}%
%\dedicatory{}%
%\commby{}%
% ----------------------------------------------------------------
\begin{abstract}
A condition on an affine central subalgebra $Z$ of a noetherian
algebra $A$ of finite Gelfand-Kirillov dimension, which we call
here \emph{unruffledness}, is shown to be equivalent in some
circumstances to the flatness of $A$ as a $Z$-module.
Unruffledness was studied by Borho and Joseph in work on
enveloping algebras of complex semisimple Lie algebras, and we
discuss applications of our result to enveloping algebras, as well
as beginning the study of this condition for more general
algebras.
\end{abstract}
\maketitle
% ----------------------------------------------------------------
\section{Introduction} \label{intro}\subsection{}\label {getgo}
Let $A$ be a noetherian algebra, finitely generated over the
uncountable algebraically closed field $k$, and let $Z$ be a
finitely generated subalgebra of the centre of $A$, such that the
nonzero elements of $Z$ are not zero divisors in $A$. The central
problem addressed in this paper is: Can we find easily checkable
conditions to ensure that $A$ is a flat $Z-$module? Of course this
question can and should be approached locally, one maximal ideal
of $Z$ at a time, but a global form of our main result
(\ref{statement}) states:

\begin{thm}
Let $A$ and $Z$ be as above, and suppose that $A$ is
Cohen-Macaulay and that $Z$ is smooth, with $\m A \neq A$ for all
maximal ideals $\m$ of $Z$. Then $A$ is a flat $Z-$module if and
only if $Z$ is unruffled in $A$.
\end{thm}

Several terms in the above statement need some explanation, which
we give in the next two paragraphs, before turning to motivation
and applications.

\subsection{}\label{cohmac} Our results and proofs are couched in the setting of
algebras of finite Gelfand-Kirillov dimension, denoted $\gk (-)$,
which we assume exists for all $A-$modules and satisfies various
standard desirable properties as listed in (\ref{stand}). The {\em
grade} $j_A(M)$ of a finitely generated $A-$module $M$ is defined
to be the least integer $j$ such that $\mathrm{Ext}^j _A (M, A)$
is non-zero, or $+\infty$ if no such integer exists; we'll simply
write $j(M)$ when the algebra $A$ is clear from the context. The
algebra $A$ is {\em Cohen-Macaulay} if
$$ \gk(A) \eq j(M) + \gk (M)
$$ for all non-zero finitely generated $A-$modules $M$. (Here and
throughout, ``module" will mean ``left module" when no other
qualification is given; so the above definition should strictly
speaking be ``left Cohen-Macaulay".) To say that $Z$ is {\em
smooth} simply means that $Z$ has finite global (homological)
dimension, or equivalently that its maximal ideal space $\Z$ is
smooth.

\subsection{} \label{unruffdef} Let $A$ and $Z$ be as in
(\ref{getgo}), and let $\mathfrak{m}$ be a maximal ideal of $Z$.
Denote the field of fractions of $Z$ by $Q(Z).$ Then $Z$ is said
to be \emph{unruffled at} $\m$ \emph{in} $A$ if
\begin{equation}\gk(A/\m A) \quad = \quad \mathrm{GK-dim}_{Q(Z)}(A \otimes_Z
Q(Z));\label{unruff} \end{equation} and $Z$ is \emph{unruffled in}
$A$ (or $A$ is \emph{unruffled over} $Z$) if (\ref{unruff}) holds
for all maximal ideals $\m$ of $Z$. The concept, although not the
name, is due to Borho and Joseph \cite[5.8]{BJ}, who showed there
that every prime factor of the enveloping algebra of a complex
semisimple Lie algebra is unruffled over its centre. Indeed,
following the suggestion of \cite[5.8]{BJ}, a secondary aim of
this paper is to begin to investigate the significance of the
unruffled hypothesis on an algebra and a central subalgebra. Our
reason for proposing the adjective ``unruffled'' is a result of
Borho \cite{Bo2}, which shows that the crucial feature of an
unruffled extension $Z \subseteq A$ is that $\gk (A/\m A)$ is
\emph{constant} as $\m$ ranges through $\Z$. Because Borho's
discussion is set in the specific context of enveloping algebras,
we shall derive a version of his result as Lemma (\ref{constant}).
To do so in the proper generality we need to recall in
(\ref{generic}) some ideas about generic ideals of algebras over
uncountable fields, which go back to work of Borho from the 1970s,
\cite[Section 4]{Bo1}. Laying out this material in a general
setting may have some independent interest.

In (\ref{ruffex}) we discuss a number of examples and non-examples
of unruffled extensions $Z \subseteq A$, and explain how our main
result collapses to a well-known theorem when $A$ is commutative.

\subsection{} \label{motivate} Pairs $Z \subseteq A$ of algebras
satisfying the hypotheses of (\ref{getgo}) arise naturally and
frequently. For example, a flat family of deformations of an
affine noetherian algebra $B$ may be exhibited as a pair $Z
\subseteq A$ of algebras as in (\ref{getgo}), with $A/\m A \cong
B$ for some particular maximal ideal $\m$ of $Z$, the deformations
of $B$ being the algebras $A/\m' A$ got by varying $\m$ across
$\Z$. Flatness of the family corresponds to $A$ being a flat
$Z-$module, so the theorem reveals that, at least in the presence
of mild hypotheses on $A$ and $Z$, this is equivalent to constancy
of the GK-dimensions of the deformed algebras.

A second major source of motivating examples is the concept of a
stratification of the prime or primitive spectrum of an algebra
$R$ into ``classically affine strata". The most clearcut examples
are given by quantum $n-$space \cite{GL}, and more generally when
$R = \mathcal{O}_q (G)$ is the quantised coordinate ring of a
semisimple group $G$ at a generic parameter $q$, \cite{HL,J1}. In
these examples the primitive spectrum $\chi$ of $R$ is the
disjoint union of finitely many locally closed subsets $\chi_w$,
and each stratum $\chi_w$ is homeomorphic to a torus. The
homeomorphism is afforded by induction $\m \mapsto \m A_w,$ where
$A_w$ is a localisation of a factor of $R$ and $\m$ is a maximal
ideal of $Z_w$, the Laurent polynomial algebra which is the centre
of $A_w.$

In a parallel mechanism, many naturally occurring algebras $R$
which are finite modules over their centres have maximal ideal
spectra which can be stratified into finitely many Azumaya strata
\cite[\S 5]{BroGor}, \cite{DCP} - including, for example,
quantised coordinate rings at a root of unity and symplectic
reflection algebras in the PI case. The point we want to make here
is not so much that the results of the present paper can
contribute anything to an understanding of Azumaya stratifications
- they can't! - but rather that some aspects of the Azumaya
stratified setting may point towards phenomena which are more
generally true. (Recall, for example, that an Azumaya algebra is
always projective over its centre.)

A third class of primitive ideal stratifications provides one of
our main motivations. Namely, let $R$ be the enveloping algebra
$U(\g)$ of a finite dimensional complex semisimple Lie algebra
with adjoint group $G$. In a series of papers Borho,
\cite{Bo2,Bo3,Bo4,Bo5}, and latterly Borho and Joseph, \cite{BJ},
have studied $\chi$, the space of primitive ideals of $R$, by
defining and studying ``generalised Dixmier maps" from subsets of
$\g^* /G$ to subsets of $\chi$. The subsets in question are the
{\em sheets} (of $\g^*/G,$ resp., of $\chi$). The most desirable
scenario - sometimes valid, sometimes not - is that a sheet
$\mathcal{S}$ in $\chi$ should (roughly speaking) consist of the
inverse images in $R$ of the ideals of a certain prime factor ring
$A$ which are generated by the maximal ideals of the centre $Z$ of
$A$. For a more detailed description of this theory and the
relevance of our results to various questions of Borho and Joseph,
see (\ref{questions}).

\subsection{} \label{apps} In Section 6 we discuss a number of
applications of the main theorem (\ref{statement}). In
(\ref{smithy}) we show that a GK-dimension inequality of Smith and
Zhang \cite{SmZh}, which is used in the proof of (\ref{statement})
and which is well-known to be strict in general, is in fact an
equality in the presence of the Cohen-Macaulay hypothesis. In
(\ref{kostant}) we derive yet another proof of the theorem of
Kostant that the enveloping algebra of a complex semisimple Lie
algebra is free over its centre, and develop this into a necessary
and sufficient criterion for arbitrary enveloping algebras. As
already mentioned, in (\ref{questions}) we explore the relevance
of (\ref{statement}) for Borho and Joseph's work on sheets of
primitive ideals. Finally, in (\ref{factors}) some preliminary
results are proved about the behaviour of the unruffled property
under factoring by a centrally generated prime ideal; and on the
way some information is produced about the how the Cohen-Macaulay
property behaves under factorisation.

\subsection{} \label{layout} As already indicated, Section 2
contains a discussion of the unruffled property and information
about ideals in general position, Section 5 contains the statement
and proof of the main theorem, and Section 6 contains
applications. The method of proof of the main theorem is
homological, and exploits a notion of depth for $Z-A-$bimodules
which are finitely generated as $A-$modules. The necessary theory
is set up in Section 3, and some technical lemmas on depth in the
presence of the unruffled hypothesis are proved in Section 4. The
final short section, Section 7, lists some questions and
suggestions for further work arising from the results described in
this paper.

\section{Unruffled extensions}\label{unruff}
\subsection{Standing hypotheses}\label{stand} We'll assume
throughout this paper that $A$ denotes an affine noetherian
algebra over the algebraically closed field $k$, and that $Z$ is
an affine subalgebra of the centre of $A$. We assume that $Z$ is a
domain whose nonzero elements are not zero divisors in $A$, as
will be the case if, for example, $A$ is prime. We write $Q(Z)$
for the field of fractions of $Z$. The maximal ideal spectrum of
$Z$ will be denoted by $\mathcal{Z}$. The Gelfand-Kirillov
dimension over $k$, denoted $\mathrm{GK-dim}_k(-)$, will be
assumed to exist for all $A$-modules, and to have the usual
desirable properties of being exact and partitive, and taking
values in the non-negative integers, as discussed in \cite{KL},
for example. Let the Gelfand-Kirillov dimensions of $A$ and $Z$ be
$n$ and $d$ respectively.

\subsection{Ideals in general position} \label{generic} As explained
in (\ref{unruffdef}), in this paragraph and the next we recall
some ideas of Borho \cite{Bo1}.

In this paragraph we assume that \begin{equation} k \textit{ is
uncountable and the }k\textit{-algebra }A\textit{ of (\ref{stand})
is finitely related.} \label{uncount} \end{equation} That is, we
assume that there exists a free $k$-algebra $F = k\langle f_1,
\ldots , f_t \rangle$ of finite rank $t$ and a finitely generated
ideal $I$ of $F$ with $F/I \cong A.$ Let $r_1, \ldots , r_m$ be a
set of generators for $I$, and for $i = 1, \ldots , m$ write
\[ r_i \eq \Sigma_{j=1}^{r(i)} \lambda_{ij}\Phi_j \]
where $\Phi_j$ are words in the free generators $f_1, \ldots ,
f_t$ of $F$ and $\lambda_{ij} \in k.$ Similarly, choose elements
$z_1, \ldots , z_r$ of $F$ whose images in $A$ generate $Z$, and
write $$ z_s \eq  \Sigma_{u=1}^{e(s)} \mu_{su}\Psi_u$$ for $s = 1,
\ldots , r$, where $\Psi_u$ are words in $f_1, \ldots , f_t$ and
$\mu_{su} \in k$. Let $k_0$ be the prime subfield of $k$ and set
$k' = k_0(\lambda_{ij}, \mu_{su} : 1 \leq i \leq m, 1 \leq j \leq
r(i), 1 \leq s \leq r, 1 \leq u \leq e(s)),$ a countable subfield
of $k$. Set $F' = k'\langle f_1, \ldots , f_t \rangle$ and $I' =
\sum_{i=1}^n F'r_i F'$, so we can define
\[ A' \quad := \quad F'/I'; \]
and set $Z'$ to be the $k'$-subalgebra of $A'$ generated by the
images in $A'$ of $z_1, \ldots , z_r.$ Thus $A = A' \otimes_{k'}
k$, so that $A'$ is a prime noetherian affine $k'$-algebra.
Clearly, $Z = Z' \otimes_{k'} k.$ (In the case where $Z = Z(A)$ we
can simply take $k' = k_0(\lambda_{ij})$ and $Z' = Z(A').$)

An ideal $\m$ of $Z$ is said to be \emph{in general position} if
$\m \cap Z' = 0.$ Versions of the following results, with similar
proofs, were obtained by Borho \cite[4.5c]{Bo1}, \cite[2.2,
2.3]{Bo2} for the case when $A$ is a prime factor of a complex
semisimple Lie algebra, and with the stronger hypothesis that
$Q(Z)A$ is a simple ring for Proposition
(\ref{stand}).3.\footnote{In fact there is a problem with part of
the argument in \cite[2.2]{Bo2}. Contrary to what is said there,
it isn't true that, for an ideal $\p$ of $Z$ in general position,
$Z' \setminus 0$ consists of regular elements $modulo(\p A)$, even
when $\p$ is semiprime, as the example $A = Z = \C [X], \quad A' =
Z' = \Q [X], \quad \p = \langle X(X - \pi ) \rangle$ makes plain.
Once this is realised, it's not hard to see that \cite[Proposition
2.2(1)]{Bo2} is false, and that the best one can say is (using
\cite[Proposition 2.2(2)]{Bo2}) that if $\p$ is semiprime with all
primes of $Z$ minimal over $\p$ in general position, then $\p A$
is semiprime.}

\begin{lem} If $\p$ is a prime ideal of $Z$ in general position
then the set \[\{ \m \in \Z : \p \subseteq \m,\quad \m \textit{ in
general position}\}\] is dense in $\mathcal{V}(\p) = \{\m \in \Z :
\p \subseteq \m \}.$
\end{lem}
\begin{proof} We may assume that $\p$ is not maximal, so that
$\mathcal{V}(\p)$ is an uncountable set. On the other hand, since
$Z'$ is countable the set
\[ \mathcal{S} \quad := \quad \cup_{z \in Z' \setminus \p}
\mathcal{V}(\p + zZ) \] is a countable union of closed proper
subsets of $\mathcal{V}(\p)$, and so does not cover
$\mathcal{V}(\p)$ by \cite[3.11]{BJa}. If $\mathcal{V}(\p)
\setminus \mathcal{S}$ were not dense then we'd have covered
$\mathcal{V}(\p)$ by a countable union of proper closed subsets,
again contradicting \cite[3.11]{BJa}. So $\mathcal{V}(\p)
\setminus \mathcal{S}$ must be dense in $\mathcal{V}(\p)$.
\end{proof}

\begin{prop} Retain the hypotheses and notation introduced in
(\ref{stand}) and in (\ref{uncount}), and let $Q(Z')$ denote the
quotient field of $Z'$.

1. $ A \quad \cong \quad Z \otimes_{Z'} A'$,\\ and hence
\begin{equation} Q(Z')A \quad \cong \quad Q(Z')Z \otimes_{Q(Z')}
Q(Z')A'. \label{gen} \end{equation} In particular, $Q(Z')A$ is a
free module over $Q(Z')Z,$ with basis including $\{1\},$ so
$Q(Z')Z$ is a direct summand of $Q(Z')A$.

2. Assume that $A$ is semiprime. If $\p$ is a prime ideal of $Z$
in general position, then $\p A$ is a semiprime ideal.

3. Assume that $A$ is prime and that $Q(Z)$ is the centre of the
simple artinian Goldie quotient ring $Q(A).$ If $\p$ is a prime
ideal of $Z$ in general position, then $\p A$ is a prime ideal.

4. Assume that $Q(Z)A$ is simple (so $A$ is prime). If $\m$ is a
maximal ideal of $Z$ in general position, then $\m A$ is a maximal
ideal.
\end{prop}

\begin{proof} 1. By the associativity of the tensor product,
\begin{equation} A \quad \cong \quad k \otimes_{k'} A' \quad \cong
\quad k \otimes_{k'} Z' \otimes_{Z'} A' \quad \cong \quad Z
\otimes_{Z'} A'. \label{assoc} \end{equation}
Localising these
isomorphisms at the central regular elements $Z' \setminus 0$ of
$A$, we find
\begin{align*} Q(Z')A \quad &= \quad  Q(Z') \otimes_{Z'} A \\& \cong \quad
Q(Z') \otimes_{Z'} (Z \otimes_{Z'} A') \\& \cong \quad Q(Z')Z
\otimes_{Z'} A' \\ &= \quad (Q(Z')Z \otimes_{Q(Z')}Q(Z'))
\otimes_{Z'} A'\\ & \cong \quad Q(Z')Z \otimes_{Q(Z')}Q(Z')A'.
\end{align*}
Since $Q(Z')A'$ is a free module over the field $Q(Z')$, the last
statement in 1. is immediate from the above isomorphisms.

2. Suppose that $A$ is semiprime, and that $\p$ is a prime ideal
of $Z$ in general position. By (\ref{gen}) and the freeness
statement in 1.,
\begin{equation} Q(Z')A\p \quad \cong \quad Q(Z')\p
\otimes_{Q(Z')}Q(Z')A', \label{pcong} \end{equation} and the
elements of $Z' \setminus 0,$ being regular $\mathit{modulo}(\p),$
are regular $\mathit{modulo}(\p A):$ for notice that, using the
last part of 1.,
\[ Q(Z')A\p \cap Z \eq Q(Z')\p \cap Z \eq \p. \]
Factoring (\ref{gen}) by (\ref{pcong}), and abusing notation
slightly by writing $Q(Z')(A/\p A)$ for the partial quotient ring
of $A/\p A$ with respect to the set $(Z' + \p A/ \p A) \setminus
0_{A/\p A},$ we obtain the isomorphism in
\begin{equation} Q(Z')(A/\p A)  \cong  Q(Z')(Z/\p
)\otimes_{Q(Z')} Q(Z')A' \subseteq  Q(Z/\p) \otimes_{Q(Z')}
Q(Z')A'. \label{blackeyed}\end{equation} The inclusion in
(\ref{blackeyed}) again follows by freeness (of $Q(Z')(A/\p A)$
over $Q(Z')(Z/\p )$), and shows that
\begin{equation} Q(Z/\p) \otimes_{Q(Z')} Q(Z')A' \textit{ is a
partial quotient ring of } A/\p A. \label{partial} \end{equation}
Since $Q(Z')A'$ is semiprime, so too is $Q(Z/\p) \otimes_{Q(Z')}
Q(Z')A'$, by \cite[3.4.2]{Dix}. But, from (\ref{blackeyed}), we
see that $Q(Z/\p) \otimes_{Q(Z')} Q(Z')A'$ is generated over $A/\p
A$ by central elements. Hence $A/\p A$ must also have no non-zero
nilpotent ideals, as required.

3. Suppose now that $A$ is prime and that $\p$ is as in 2. Suppose
that $Q(Z)$ is the centre of $Q(A)$. Then $Q(Z')$ is the centre of
$Q(A')$, since, by (\ref{gen}),
\[ Q(A) \quad \cong \quad Q(Z) \otimes_{Q(Z')} Q(A'). \]
So by \cite[proof of 7.3.9]{Pas} $Q(Z/\p ) \otimes_{Q(Z')} Q(A')$
is simple, and hence, being noetherian, it has a simple artinian
quotient ring by Goldie's theorem \cite{MR}. That $\p A$ is prime
now follows from (\ref{partial}).

4. Suppose now that $\m$ is a maximal ideal of $Z$ in general
position. Since $\m \cap Z' = 0,$ the map $Z \To Z/\m$ induces a
homomorphism from $Q(Z')Z$ to $k$, so $Q(Z') \subseteq k$ and we
can form the tensor product $k \otimes_{Q(Z')} Q(Z')A'.$ Thus
(\ref{blackeyed}) simplifies to
\begin{equation} A/\m A \eq Q(Z')(A/\m A) \quad \cong \quad k \otimes_{Q(Z')}
Q(Z')A'. \label{peas} \end{equation}  Suppose now that $Q(Z)A$ is
simple. From (\ref{assoc}), $Q(Z)A \cong Q(Z) \otimes_{Q(Z')}
Q(Z')A',$ so that $Q(Z')A'$ is also simple. Simplicity of $A/\m A$
follows from this by (\ref{peas}) and \cite[proof of 7.3.9]{Pas}.
\end{proof}

\subsection{Generic constancy of GK-dimension}\label{constant} As
already explained, the following result was obtained by Borho and
Joseph for factors of enveloping algebras, with the same proof. It
seems reasonable to suspect the truth of a stronger result -
namely, that the set of unruffled maximal ideals of $Z$ in $A$
contains a non-empty Zariski-open subset of $\Z$.

\begin{lem} \cite[5.8]{BJ} Keep the hypotheses on $Z$ and $A$ from
(\ref{stand}) and (\ref{generic}), (but there is no need to assume
that $Q(Z)A$ is simple). Let $\m$ be a maximal ideal of $Z$ in
general position. Then
\[ \gk (A/\m A) \eq \mathrm{GK-dim}_{Q(Z)}(Q(Z) \otimes_Z A). \]
\end{lem}
\begin{proof} Associativity of the tensor product yields
\begin{align} Q(Z) \otimes_{Q(Z')} Q(Z')A' \quad &= \quad Q(Z) \otimes_{Q(Z')}
Q(Z')\otimes_{Z'} A' \nonumber \\
&=\quad (Q(Z) \otimes_Z Z) \otimes_{Z'} A' \eq Q(Z) \otimes_Z A.
\label{ass} \end{align} From (\ref{peas}) we get
\begin{align*} \gk (A/\m A) \quad &= \quad
\mathrm{GK-dim}_{Q(Z')}(Q(Z')A') \\
&= \quad \mathrm{GK-dim}_{Q(Z)}(Q(Z) \otimes_{Q(Z')} Q(Z')A') \\
&= \quad \mathrm{GK-dim}_{Q(Z)}(Q(Z) \otimes_Z A),
\end{align*}
where the final equality is given by (\ref{ass}).
\end{proof}

\subsection{Unruffled and ruffled examples} \label{ruffex} Recall
that, where nothing is said to the contrary, hypotheses
(\ref{stand}) are assumed to hold throughout.

\textbf{1.} \textit{The case where $A$ is a finitely generated
$Z$-module.} It's clear that if $A$ is a finitely generated
$Z$-module, then $A$ is unruffled over $Z$. One only needs to note
that if $\m$ is a maximal ideal of $Z$ then $\m A$ is a proper
ideal, which can be seen by inverting the regular elements $Z
\setminus \m$ in $A$ and appealing to Nakayama's lemma.

\textbf{2.} \textit{Prime factors of semisimple enveloping
algebras are unruffled over their centres.} If $A =
U(\mathfrak{g})/P$ is a prime factor of the enveloping algebra of
a finite dimensional complex semisimple Lie algebra
$\mathfrak{g}$, then $A$ is an unruffled extension of its centre
by \cite[Corollary 5.8]{BJ}. The existing proof of this fact is
rather deep, depending as it does on the description of $P$ as
induced from a rigid primitive ideal of the enveloping algebra of
a Levi subalgebra $\mathfrak{l}$ of a parabolic subalgebra of
$\mathfrak{g}$ combined with an irreducible subset of the centre
of $\mathfrak{l}$.

\textbf{3.} \textit{The commutative case.} Suppose that all the
assumptions of (\ref{stand}) hold, but in addition $A$ is
commutative and Cohen-Macaulay, so $Z$ is now an arbitrary affine
subalgebra of $A$. Routine local-global yoga applied to
\cite[Theorem 18.16b and Corollary 13.5]{E} easily yields our main
result in this commutative setting: \emph{If $\m$ is a smooth
point of $\Z$, then $Z$ is unruffled in $A$ at $\m$ if and only if
$\m A \neq A$ and $A_{\m} := A \otimes_Z Z_{\m}$ is a flat
$Z_{\m}$-module.}

An instructive example to consider here is the subalgebra $Z = \C
[x,xy]$ of the commutative polynomial algebra $A = \C [x,y].$ One
easily confirms that, for a maximal ideal $\m$ of $Z$, $A_{\m}$ is
a flat $Z_{\m}$-module if and only if $\m \neq \langle x, xy
\rangle,$ while $\m$ is unruffled in $A$ if and only if $\m \neq
\langle x, xy - \lambda \rangle,$ for $\lambda \in \C.$

If one assumes, in addition to the commutativity of $A$, that $A$
is a finitely generated $Z$-module, then, noting (\ref{ruffex}).1,
one recovers from Theorem (\ref{getgo}) the familiar fact
\cite[Corollary 18.17]{E} that a commutative affine Cohen-Macaulay
domain is projective over any smooth subring over which it's a
finitely generated module.

\textbf{4.} \textit{Enveloping algebras of solvable Lie algebras
are not always unruffled over their centres.} Let $\mathfrak{g}$
be the complex solvable Lie algebra with basis $x,y,z,t$, such
that
$$[t,x] = x, \quad [t,y] = -y, \quad [t,z] = -z,$$ and all other
brackets are 0. Let $A = U(\mathfrak{g})$ and let $Z$ be the
centre of $A$. Thus $A = R[t; \delta]$ where $R = \C[x,y,z]$ is a
commutative polynomial algebra and $\delta$ is a derivation. One
calculates easily that $Z$ is contained in $R$, so that $Z$
consists of the $\delta$-invariants in $R$. Since $\delta$ acts
semisimply on $R$, with the eigenvector $x^iy^jz^{\ell}$ having
eigenvalue $i - j -\ell$, it follows that $$Z = \{\Sigma \C
x^iy^jz^{\ell}: i = j +\ell \} = \C [xz,xy], $$ a polynomial
algebra in two variables. For $a,b,\in \C$, let  $\m_{a,b}$ denote
the maximal ideal $\langle xy - a, xz - b \rangle$ of $Z$. It's
routine to check that
\[ \mathrm{GK-dim}_{\C}(A/\m_{a,b}) \eq 2 \]
for $(a,b) \neq (0,0)$, while $A/\m_{0,0}$ maps onto $\C[y,z][t;
\delta ]$, so that
\[ \mathrm{GK-dim}_{\C}(A/\m_{0,0}) \eq 3.\]
Thus $Z$ is not unruffled in $A$ at $\m_{0,0}$.

\textbf{5.} \textit{Left noetherian PI-rings are not always
unruffled over their centres.} Let $t$ and $s$ be indeterminates,
and define
\[ A \quad = \quad \left[
         \begin{array}{rr}
              k[t,t^{-1},s] & k[t,t^{-1},s] \\
              0 &  k[t]
          \end{array} \right], \]
where $k[t,t^{-1},s]$ is a right $k[t]$-module via the embedding
of the second algebra in the first. Thus $A$ is a left noetherian
affine PI algebra, but is not semiprime and is not right
noetherian. Set $Z$ to be the centre of $A$, which is easily
checked to be the set of scalar matrices and so isomorphic to
$k[t]$. Thus $Z \setminus 0$ consists of regular elements of $A$,
and $Q(Z) \cong k(t),$ with
\[ A \otimes_Z Q(Z) \quad = \quad \left[
         \begin{array}{rr}
              k(t)[s] & k(t)[s] \\
              0 &  k(t)
          \end{array} \right].\]
Thus
\[ \mathrm{GK-dim}_{Q(Z)}(A \otimes_Z Q(Z) ) \eq 1. \]
Consider the maximal ideal
\[ \m \quad = \quad \left[
         \begin{array}{rr}
              t & 0 \\
              0 &  t
          \end{array} \right]Z \]
of $Z$. One easily calculates that $A/\m A \cong k$, so that
\[\gk (A/\m A) = 0. \]
So $Z$ is not unruffled in $A$ at $\m$.

\subsection{Inequalities for unruffled extensions}\label{ruffineq}
In the presence of flatness the following lemma shows that the
strict inequality of GK-dimensions in Example (\ref{ruffex}).5 is
the \emph{only} direction in which unruffledness can fail. But
Example (\ref{ruffex}).3 (with $\m = \langle x,xy \rangle$) shows
that the flatness hypothesis in the lemma is needed.

\begin{lem} Let $A$ and $Z$ be as in (\ref{stand}), and let $\m
\in \mathcal{Z}$. Suppose that $A_{\m}$ is a flat $Z_{\m}-$module.
Then
\begin{equation}
\gk (A/\m A) \quad \leq \quad \mathrm{GK-dim}_{Q(Z)}(A \otimes_Z
Q(Z)). \label{ruff1}
\end{equation}
\end{lem}
\begin{proof} Denote $Q(Z)$ by $Q$. Let $V$ be a finite dimensional $k$-vector space
which generates $A$ as a $k$-algebra. So the image $\overline{V}$
of $V$ in $A/\m A$ [resp. the image of $V$ in $A \otimes_Z Q$]
generates $A/\m A$ [resp. $A \otimes_Z Q$] as a $k-$ [resp. $Q-$]
algebra. It will therefore be enough to show that, for all $i \geq
1$,
\begin{equation}
\mathrm{dim}_k (\overline{V}^i) \quad \leq \quad \mathrm{dim}_Q
(QV^i).\label{spencer}
\end{equation}
Suppose then that $u_1, \ldots , u_t$ are elements of $V^i$ such
that $\Sigma_{j=1}^t q_j u_j = 0$, where $q_j \in Q$, not all
zero. We claim that $\overline{u}_1, \ldots , \overline{u}_t$ are
$k$-linearly dependent elements of $\overline{V}^i$. It's clear
that this will prove (\ref{spencer}).

Multiplying by a suitable element of $Z$ and discarding those
$u_j$ for which $q_j = 0$, we get
\[ \Sigma_{j=1}^t z_j u_j \quad = \quad 0, \]
with each $z_j$ a non-zero element of $Z$. Fix a maximal ideal
$\m$ of $Z$. Choose $\ell \geq 1$, $\ell$ minimal such that there
exists $j$ with $z_j \notin \m^{\ell}$. (Note that $\ell$ exists
by the Krull Intersection Theorem \cite[Corollary 5.4]{E}, $Z$
being a noetherian domain.) Thus
\begin{equation}
\Sigma_{j=1}^t (z_j + \m^{\ell})(u_j + \m^{\ell} A) \quad = \quad
0 \label{blues}
\end{equation}
in $A/ \m^{\ell} A$, with not all the $z_j + \m^{\ell}$ equal to
0. Let $\gamma_1, \ldots , \gamma_p$ be a $k$-basis for $\m^{\ell
- 1}/\m^{\ell}$, and write $z_j + \m^{\ell} = \Sigma_{r=1}^p
\gamma_r \lambda_{jr}$, for $\lambda_{jr} \in k$. Thus
(\ref{blues}) gives
\[ \Sigma_j (\Sigma_r \gamma_r \lambda_{jr} + \m^{\ell})(u_j +
\m^{\ell}A) \quad = \quad 0. \] That is,
\begin{equation}
\Sigma_r\gamma_r (\Sigma_j \lambda_{jr} + \m^{\ell})(u_j +
\m^{\ell}A) \quad = \quad 0. \label{jon}
\end{equation}
Now the linear independence of $\{\gamma_r\}$ in $\m^{\ell -
1}/\m^{\ell}$ over $k$ implies, thanks to the flatness hypothesis,
linear independence of $\{\gamma_r \}$ in $\m^{\ell -
1}A/\m^{\ell}A$ over $A/\m A$. Hence (\ref{jon}) shows that, for
each $r$,
\[ \Sigma_{j=1}^t \lambda_{jr} u_j \quad \in \quad \m A.\]
For some $r$, there exists $j$ with $\lambda_{jr} \neq 0$. So the
result is proved.
\end{proof}

\section{Homological equipment}\label{hom}
\subsection{Depth} \label{depth} We need to extend the standard
notion of depth from commutative algebra. The classical definition
(as in for example \cite[page 425]{E}) begins with a commutative
noetherian ring $R$, an ideal $I$ of $R$, and a finitely generated
$R$-module $M$ with $MI \neq M$, and defines the \emph{depth} of
$I$ on $M$ to be the length of a maximal $M$-sequence of elements
of $I$. (Recall that an $M$-\emph{sequence} is a sequence $\{x_1,
\dots , x_n \}$ of elements of $R$ such that $x_i$ is not a zero
divisor on $M/\Sigma_{j=1}^{i-1} x_j M$, for $i = 1, \dots , n$;
the \emph{length} of the $M$-sequence is then $n$.) Crucial to the
usefulness of this definition is \cite[Theorem 17.4]{E}, which
guarantees that any two such maximal $M$-sequences have the same
length, and that this number can be read off from an appropriate
Koszul complex.

We extend the above definition by allowing the $R$-module $M$ to
be not necessarily finitely generated, but we still insist that
$MI \neq M$, and we require $M$ to be an $S-R$-bimodule with $S$ a
left noetherian ring and $M$ a finitely generated $S$-module. With
this definition, the analogue of \cite[Theorem 17.4]{E}, which we
state below and prove in (\ref{depthproof}), remains true. Write
$R^{(n)}$ for the direct sum of $n$ copies of $R$. For elements
$x_1, \ldots , x_n$ of the commutative noetherian ring $R$, we
denote by $K_R (x_1, \ldots , x_n)$, or by $K (x_1, \ldots , x_n)$
when the ring is clear from the context, the Koszul complex
\[ 0 \To R \To R^{(n)} \To \wedge^2 R^{(n)} \To \cdots \To
\wedge^{i} R^{(n)} \overset{d_{\mathbf{x}}}{\To}
\wedge^{i+1}R^{(n)} \To \cdots \wedge^n R^{(n)} \To 0, \] with
$\mathbf{x} = (x_1, \ldots , x_n) \in R^{(n)}$ and
$d_{\mathbf{x}}(a) = \mathbf{x} \wedge a$, \cite[pages 423-4]{E}.

\begin{thm} Let $R$, $I$, $S$ and $M$ be as stated above, and
suppose that $I = \Sigma_{i=1}^n x_i R$. Let $r$ be a non-negative
integer. If
\[ H^j (M \otimes_R K(x_1,
\ldots , x_n)) \quad = \quad 0 \] for $j < r$, while
\[ H^r (M \otimes_R K(x_1,
\ldots , x_n)) \quad \neq \quad 0 ,\] then every maximal
$M$-sequence in $I$ has length $r$.
\end{thm}

\subsection{Lemma:}\label{primes}\emph{Let} $R$, $S$ \emph{and} $M$ \emph{be as stated in}
(\ref{depth}), \emph{with} $M \neq 0$.

\emph{1. The set of zero divisors of} $R$ \emph{on} $M$ \emph{is
equal to the union of a finite set of prime ideals of} $R$.

\emph{2. If} $I$ \emph{is an ideal of} $R$ \emph{which consists of
zero divisors on }$M$\emph{, then there exists a prime ideal
}$\mathfrak{p}$ \emph{of }$R$ \emph{with }$I \subseteq
\mathfrak{p}$\emph{, and} $0 \neq m \in M$\emph{, with
}$m\mathfrak{p} = 0.$
%\end{lem}
\begin{proof} Since $M$ is an $S-R$-bimodule and is left
noetherian, it has by \cite[Proposition 4.4.9]{MR} an affiliated
series of prime ideals $\{\mathfrak{p_1}, \ldots ,
\mathfrak{p_m}\}$ as an $R$-module, in the sense of
\cite[4.4.6]{MR}, such that no element of $R \setminus
(\cup_{i=1}^m \mathfrak{p_i})$ is a zero divisor on $M$. Thus
\[ I \quad \subseteq \quad \cup_{i=1}^m \mathfrak{p_i}. \]
By the prime avoidance property \cite[Lemma 3.3]{E} there exists
$j$, $1 \leq j \leq m$, such that $I \subseteq \mathfrak{p_j}$.
Since $\{m \in M : m\mathfrak{p_j} = 0 \}$ is a non-zero submodule
of $M$, the lemma is proved.
\end{proof}

\subsection{}\label{koszul} For the most part the proof of Theorem (\ref{depth})
follows the classical approach as in \cite[Section 17.3, proof of
17.4 ]{E}. Thus \cite[Proposition 17.9 and Corollaries 17.10.
17.11]{E} don't involve the module $M$ and so apply unchanged
here. But we require an improved version of \cite[Corollary
17.12]{E}.
\begin{prop} Let $R$, $I$, $S$ and $M$ be as in (\ref{depth}),
with $I = \Sigma_{i=1}^n x_i R$. Suppose that $r$ is a
non-negative integer and that $x_1, \ldots , x_r$ is an
$M$-sequence. Then \begin{equation} H^r(M \otimes_R K(x_1, \ldots
, x_n)) \quad = \quad \{m \in M : mI \subseteq \Sigma_{i=1}^r Mx_i
 \} / \Sigma_{i=1}^r Mx_i. \label{Kos1} \end{equation} Hence,
for $j < r$,
\begin{equation} H^j(M \otimes_R K(x_1, \ldots , x_n)) \quad =
\quad 0, \label{Kos2} \end{equation}  and if $\{x_1, \ldots , x_r
\}$ is a maximal $M$-sequence in $I$ then \begin{equation} H^r(M
\otimes_R K(x_1, \ldots , x_n)) \quad \neq \quad 0. \label{Kos3}
\end{equation}
\end{prop}

\begin{proof} We prove (\ref{Kos1}) by induction on $r$; for $r =
0$, the statement follows from the definition of the Koszul
complex. Now suppose that $r > 0$, with the result proved for
smaller values of $r$. We use here induction on $n$, starting from
$n = r$. In this starting case, (\ref{Kos1}) states that $H^r(M
\otimes_R K(x_1, \ldots , x_r)) = M/MI$, which is clear from the
definition of the Koszul complex.

Suppose now that $n > r$, and the result is known for this $r$ and
smaller values of $n$. By the induction on $r$ we have
\[H^{r-1}(M \otimes_R K(x_1, \ldots , x_n)) = \{m \in M : mI
\subseteq \Sigma_{i=1}^{r-1}Mx_i \}/\Sigma_{i=1}^{r-1}Mx_i = 0, \]
since $x_r$ is not a zero divisor on $M/\Sigma_{i=1}^{r-1}Mx_i$.
Thus the exact sequence of \cite[Corollary 17.11]{E} yields
\begin{align}
H^{r}(M \otimes_R &K(x_1, \ldots , x_n)) \eq \nonumber \\
&\mathrm{ker}(H^{r}(M \otimes_R K(x_1, \ldots , x_{n-1}))
\overset{x_n \times}{\To} H^{r}(M \otimes_R K(x_1, \ldots ,
x_{n-1})). \label{Kos4}
\end{align} Write $N = \{m \in M : m(\Sigma_{i=1}^{n-1}x_i R)
\subseteq \Sigma_{i=1}^r Mx_i \}$. Then\begin{align} \{m \in M :
mI \subseteq \Sigma_{i=1}^r Mx_i \}&/\Sigma_{i=1}^r Mx_i \eq
\nonumber \\ &\mathrm{ker}((N/\Sigma_{i=1}^r Mx_i) \overset{x_n
\times}{ \To} (N/\Sigma_{i=1}^r Mx_i)). \label{Kos5} \end{align}
Comparing (\ref{Kos4}) with (\ref{Kos5}) proves the induction step
for (\ref{Kos1}).

Since $x_{j+1}$ is not a zero divisor on $M/\Sigma_{i=1}^j Mx_i$
for $j < r$, (\ref{Kos2}) follows at once from (\ref{Kos1}). To
prove (\ref{Kos3}), suppose that $\{x_1, \ldots , x_r \}$ is a
maximal $M$-sequence in $I$. Then $I$ is contained in the set of
zero divisors on $M/\Sigma_{i=1}^rMx_i$. By Lemma (\ref{primes})
there exists $m \in M$, $m \notin \Sigma_{i=1}^rMx_i$, such that
$mI \subseteq \Sigma_{i=1}^rMx_i$. So (\ref{Kos3}) follows from
this and (\ref{Kos1}).
\end{proof}

\subsection{Proof of Theorem (\ref{depth})}\label{depthproof} Let
$y_1, \ldots , y_s$ be a maximal $M$-sequence in $I$. By
hypothesis, $r$ is the least integer $j$ such that $H^{j}(M
\otimes_R K(x_1, \ldots , x_n)) \neq 0.$ Now $r$ is also the least
integer $j$ for which $H^{j}(M \otimes_R K(x_1, \ldots , x_n,y_1,
\ldots , y_s)) \neq 0,$ by \cite[Corollary 17.10]{E}. Since $MI
\neq M$ by hypothesis, Proposition (\ref{koszul}) shows that
$s=r,$ proving the theorem.     $\square$

\subsection{Definition of depth}\label{depdef} Let $R,S,I$ and $M$
be as in (\ref{depth}). Define the \emph{depth} of $I$ on $M$,
denoted $\mathrm{depth}(I,M)$, to be the length of a maximal
$M$-sequence in $I$. Theorem (\ref{depth}) shows that this
definition makes sense.

\subsection{Grade versus depth}\label{extdepth} We need a
noncommutative variant of one of the standard commutative
characterisations of depth, as given in \cite[Proposition
18.4]{E}, for example.

\begin{prop} Let $S$ be a noetherian ring with a noetherian
central subring $R$ and let $J$ be an ideal of $S$ with $J = (J
\cap R)S$. Let $N$ be an $S-S$-bimodule, finitely generated on
each side, with $R$ acting centrally on $N$ and $NJ \neq N$. Then
the depth of $J \cap R$ on $N$ is equal to the least non-negative
integer $r$ such that $\mathrm{Ext}^r_S (S/J, N) \neq 0.$
\end{prop}

\begin{proof} Assume that $S,R,N$ and $J$ are as stated. We prove
the theorem by induction on $\mathrm{depth}(J \cap R, N) = t$.
Suppose first that $t = 0.$ Then it's immediate from Lemma
(\ref{primes}) that $\mathrm{Hom}_S (S/J, N) \neq 0,$ as required.

Now assume that $t \geq 1$ and that the result is proved for
smaller values of the depth. Let $x \in J \cap R$ be a regular
element on $N$. We have $J(N/xN) \neq N/xN,$ and $\mathrm{depth}(J
\cap R, N/xN) = t-1$ by Theorem (\ref{depth}). So, by induction,
$\mathrm{Ext}^{t-1}_S (S/J, N/xN) \neq 0,$ but $\mathrm{Ext}^{i}_S
(S/J, N/xN) = 0$ for all $i < t-1.$ Applying $\mathrm{Hom}_S(S/J,
- )$ to the exact sequence
\[ 0 \To N \overset{x \times}{\To} N \To N/xN \To 0, \] we get for each $j \geq 1$ the exact sequence
\[ 0 \To \mathrm{Ext}^{j-1}_S (S/J, N) \To \mathrm{Ext}^{j-1}_S (S/J,
N/xN)\To \mathrm{Ext}^{j}_S (S/J, N) \To 0, \] where the  first
and last terms are $0$ because $x\mathrm{Ext}^{i}_S (S/J, N) = 0$
for all $i$. Hence we deduce that $\mathrm{Ext}^{i}_S (S/J, N) =
0$ for $i < t,$ and $\mathrm{Ext}^{t}_S (S/J, N) \neq 0,$ as
required.
\end{proof}

\subsection{Measuring the flat dimension:}\label{torlem}The following
result is standard and easy for finitely generated modules over a
commutative noetherian ring \cite[Theorem 6.8]{E}, but is false
for infinitely generated modules without some additional
hypothesis, as can be seen by taking $Z$ to be a polynomial ring
in 2 variables and $M$ to be the field of fractions of the factor
by a height one prime.

\begin{lem} Let $A$ and $Z$ be as in (\ref{stand}) and suppose that $k$ is an uncountable field. Let $M$ be a
finitely generated $A$-module which has finite flat dimension $t$
as a $Z$-module. Then
\begin{align*} t \quad &= \quad \mathrm{max} \{ r : \mathrm{Tor}^r_Z(V,M) \neq
0, V \: \mathit{a} \:
Z-\mathit{module,}\; \mathrm{dim}_k (V) < \infty \}\\
&= \quad \mathrm{max} \{ r : \mathrm{Tor}^r_Z(Z/\mathfrak{m},M)
\neq 0,
\mathfrak{m} \in \mathcal{Z} \} \\
&= \quad \mathrm{max} \{  :
\mathrm{Tor}^r_{Z_{\mathfrak{m}}}(Z/\mathfrak{m},M_{\mathfrak{m}})
\neq 0, \mathfrak{m} \in \mathcal{Z} \}.
\end{align*}
\end{lem}

\begin{proof} The second equality is an easy consequence of the
long exact sequence of $\mathrm{Tor}$, and the third is clear
since $\mathfrak{m}\mathrm{Tor}^r_Z(Z/\mathfrak{m},M) = 0.$ Since
$t$ is finite by hypothesis, it is an upper bound for the right
side of the first equality. Moreover the long exact sequence of
$\mathrm{Tor}$ also shows easily that there exists a prime ideal
$\mathfrak{p}$ of $Z$ with $\mathrm{Tor}^t_Z(Z/\mathfrak{p},M)
\neq 0.$ Choose $\mathfrak{p}$ to be maximal among such primes,
and suppose for a contradiction that $\mathfrak{p}$ is not a
maximal ideal. Let $y \in Z \setminus \mathfrak{p}$, with $y +
\mathfrak{p}$ not a unit. The exact sequence
\[ 0 \To Z/\mathfrak{p} \overset{y \times}{\To} Z/\mathfrak{p} \To Z/ \mathfrak{p}+ yZ
\To 0 \] yields
\[ \mathrm{Tor}^{t+1}_Z(Z/\mathfrak{p}+yZ,M) \To
\mathrm{Tor}^t_Z(Z/\mathfrak{p},M)\overset{y \times}{
\To}\mathrm{Tor}^t_Z(Z/\mathfrak{p},M)\To
\mathrm{Tor}^t_Z(Z/\mathfrak{p}+yZ,M),\]in which the two outer
terms are zero by our hypotheses on $t$ and $\mathfrak{p}$. Thus
multiplication by $y$ is a bijection on
$\mathrm{Tor}^t_Z(Z/\mathfrak{p},M);$ in other words,
$\mathrm{Tor}^t_Z(Z/\mathfrak{p},M)$ is a vector space over the
quotient field $Q(Z/\mathfrak{p})$ of $Z/\mathfrak{p}$. But since
$k$ is uncountable and $\mathfrak{p}$ is not maximal,
$\mathrm{dim}_k (Q(Z/\mathfrak{p}))$ is uncountable. Hence
$\mathrm{dim}_k(\mathrm{Tor}^t_Z(Z/\mathfrak{p},M))$ is also
uncountable. This, however, is impossible, since
$\mathrm{Tor}^t_Z(Z/\mathfrak{p},M)$ is a finitely generated
module over the countable dimensional $k$-algebra $A$.
\end{proof}

\section{Unruffled technicalities} \label{unruff2}
\subsection{}\label{ruffvalue} We shall assume throughout
Section \ref{unruff2} that $A$ and $Z$ satisfy the hypotheses of
(\ref{stand}), (so in particular they have GK-dimensions $n$ and
$d$ respectively). Recall that the definitions of the
Cohen-Macaulay property and of the grade $j(M)$ of an $A$-module
$M$ are given in (\ref{cohmac}).

\begin{lem} Let $A$ and $Z$ be as in (\ref{stand}) and assume that
$A$ is Cohen-Macaulay. Let $\mathfrak{m}$ be a smooth point of
$\Z$, and suppose that $\m A \neq A$ and that \begin{equation} \gk
(A/\m A) \leq \mathrm{GK-dim}_{Q(Z)}(A \otimes_Z Q(Z)).
\label{weakunruff}
\end{equation} Then \begin{equation}
\mathrm{GK-dim}_k(A/\mathfrak{m}A) \quad = \quad n-d,\label{conc}
\end{equation}
and $Z$ is unruffled in $A$ at $\m$.
\end{lem}

\begin{proof} The Cohen-Macaulay property of $A$ implies that
\begin{equation} n \quad = \quad  \mathrm{GK-dim}_k(A/\mathfrak{m}A)
+ j(A/\mathfrak{m}A), \label{CM}
\end{equation}
and we note that the validity of (\ref{CM}) is unaffected by
inverting the powers of any element of $Z \setminus \mathfrak{m}$
in $A$, by \cite[Proposition 4.2]{KL} and the fact that
$\mathrm{Ext}^j_A(A/\mathfrak{m}A, A)$ is annihilated by
$\mathfrak{m}$ for all $j$. Similarly, our desired conclusion
(\ref{conc}) is clearly unaffected by such a localisation. So we
invert in $A$ the powers of an element $x$ of $Z \setminus
\mathfrak{m}$, chosen so that in the localised ring $Z[x^{-1}]$,
$\mathfrak{m}$ is generated by a regular sequence $\{x_1, \ldots ,
x_d \}$.

By \cite[Corollary 2]{SmZh},
\begin{equation} n \quad = \quad \mathrm{GK-dim}_k(A) \quad \geq \quad
\mathrm{GK-dim}_{Q(Z)}(A\otimes_Z Q(Z)) + d. \label{smzh}
\end{equation}
By (\ref{weakunruff}) and (\ref{smzh}),
\begin{equation} n - \mathrm{GK-dim}_k(A/\mathfrak{m}A) \quad \geq
\quad d. \label{gkineq}
\end{equation}
By (\ref{CM}) and (\ref{gkineq}),
\begin{equation} j(A/\mathfrak{m}A) \quad \geq \quad d.
\label{gradeineq} \end{equation} In view of Proposition
(\ref{extdepth}) we can rewrite (\ref{gradeineq}) as
\begin{equation} \mathrm{depth}(\mathfrak{m},A) \quad \geq \quad
d. \label{depthineq} \end{equation} On the other hand the Koszul
complex $K_Z (x_1, \ldots , x_d )$ gives a $Z$-free resolution of
$Z/\mathfrak{m},$ and applying $- \otimes_Z A$ to this we see that
\[ H^d(K_Z (x_1, \ldots , x_d) \otimes_Z A) \quad = \quad
A/\mathfrak{m}A \quad \neq \quad 0. \] Thus Theorem (\ref{depth})
implies that \begin{equation} \mathrm{depth}(\m,A) \leq
d.\label{depth2ineq}\end{equation} From (\ref{depthineq}) and
(\ref{depth2ineq}), and Proposition (\ref{extdepth}) we find that
equality holds in (\ref{gradeineq}); that is, $j(A/\m A) =  d,$
and substituting this value in (\ref{CM}) gives (\ref{conc}).

Moreover, substituting (\ref{conc}) in (\ref{smzh}) yields
\begin{equation}\gk (A/\m A) \geq
\mathrm{GK-dim}_{Q(Z)}(A \otimes_Z Q(Z)), \label{weakunruff2}
\end{equation}
so that, given (\ref{weakunruff}), $Z$ is unruffled in $A$ at
$\m$.

\end{proof}

\subsection{Lemma}\label{primeintersect} \emph{Let} $Z$ \emph{and} $A$ \emph{be as in
(\ref{stand}), and suppose that }$Z$ \emph{is unruffled in} $A$.
\emph{For every prime} $\mathfrak{p}$ \emph{of }$Z$, $\p A \cap Z
= \p.$

\begin{proof} The unruffled hypothesis forces $\m A \cap Z = \m$
for every maximal ideal $\m$ of $Z$. If $\p$ is a prime ideal of
$Z$ then
\begin{align*} \p \quad \subseteq \quad \p A \cap Z \quad &\subseteq \quad \cap\{\m A \cap
Z : \p \subseteq \m \in \mathcal{Z} \}\\
&= \quad \cap\{\m  : \p \subseteq \m \in \mathcal{Z} \}\\
&= \quad \p,
\end{align*}
the last equality holding since $Z$ is affine over $k$
\cite[Theorem 4.19]{E}.
\end{proof}

\subsection{}\label{primeunruff} The next result extends one direction of the equality in Lemma (\ref{ruffvalue})
from maximal to prime ideals of $Z.$ We'll improve both
inequalities below to equalities in Theorem (\ref{factors}),
provided $Z$ is smooth and (\ref{uncount}) holds.

\begin{lem}Let
$A$ and $Z$ be as in (\ref{stand}), and suppose that $A$ is
Cohen-Macaulay. Let $\p$ be a prime ideal of $Z$ of height $\ell$
which is not in the singular locus, and suppose that $Z$ is
unruffled in $A$ at the smooth points of $\mathcal{Z}$. Then
\begin{equation} \gk (A/ \p A) \geq n - \ell \label{finis} \end{equation}
and
\begin{equation} j(A/\p A) \leq \ell. \label{gravadlax} \end{equation}
\end{lem}
\begin{proof} Since $A$ is Cohen-Macaulay of Gelfand-Kirillov dimension $n$, (\ref{finis}) and (\ref{gravadlax})
are equivalent; we prove (\ref{finis}). Suppose we invert in $A/\p
A$ the powers of an element $z$ of $Z \setminus \p$; if $z$ is not
a zero divisor {\em modulo}$(\p A)$ then $\gk (A/ \p A)$ is
unchanged by this localisation, \cite[Proposition 4.2]{KL}, while
if $z$ {\em is} a zero divisor then $\gk (A/ \p A)$ may decrease
when we invert $z$. So in proving (\ref{finis}) we may invert a
suitable element of the ideal defining the singular locus and so
arrange that $Z$ is smooth. We argue by induction on
\begin{equation}\label{deft} t \quad := \quad \gk (Z/\p) \eq d - \ell. \end{equation}
 The
starting point $t = 0$ is given by Lemma (\ref{ruffvalue}).

Suppose that $t$ is greater than 0, and that we have shown that
\begin{equation} \gk (A/ \mathfrak{q}A) \quad \geq \quad n - (\ell
+ 1) \label{qgk} \end{equation} for all primes $\mathfrak{q}$ of
height $(\ell + 1)$. We apply Lemma (\ref{primes})(1) with $M = A/
\p A$, which is a non-zero module by Lemma (\ref{primeintersect}).
The same lemma in fact tells us that $\mathrm{Ann}_Z(M) = \p,$ and
since $Z/\p$ is an affine $k$-algebra of infinite $k$-dimension,
Lemma (\ref{primes})(1) ensures that there exists $x \in Z$ with
$x + \p$ a non-unit of $Z/\p$ such that $x + \p A$ is not a zero
divisor in $A/\p A$. So by \cite[Proposition 5.1(e)]{KL},
\begin{equation} \gk (A/ \p A + xA) \quad < \quad \gk (A/\p A).
\label{gkdown} \end{equation} But $\p A + xA = (\p + xZ)A,$ and
$\gk (Z/\p + xZ) = t-1$ by the Principal Ideal Theorem
\cite[Theorem 10.1 and Corollary 13.4]{E}. Thus the induction
hypothesis (\ref{qgk}) coupled with (\ref{gkdown}) yields
(\ref{finis}). This proves the induction step and hence the lemma.
\end{proof}

\subsection{Example.} \label{Heis} \emph{Lemmas} (\ref{ruffvalue}) \emph{and}
(\ref{primeunruff}) \emph{are in general false if }$A$ \emph{is
not Cohen-Macaulay.} Consider the Heisenberg group $H$ on 2
generators,
$$H = \langle x,y,: [[x,y],x] =
[[x,y],y] = 1\rangle.$$  Set $z = [x,y]$ and let $Z$ be the
subalgebra $k\langle z \rangle$ of the group algebra $A = kH$.
Thus $ Z$ is the centre of $A$ and clearly $A$ is a free
$Z$-module. By \cite[Example 11.10]{KL}
$$ \gk (kH) \eq 4.$$ One can easily see that, for all maximal ideals $\m$ of $Z$,
\[ \gk (A/\m A) \eq 2 \eq \mathrm{GK-dim}_{Q(Z)}(A \otimes_Z
Q(Z)). \] Thus Lemma (\ref{ruffvalue}) fails here; clearly $A$ is
not Cohen-Macaulay, since, for all maximal ideals $\m$ of $Z$,
$j(A/\m A) = 1.$

\section{The main theorem}\label{Main}
\subsection{}\label{statement} After stating the result we shall
prove the first part in (\ref{proof1}) and the second in
(\ref{proof2}). Clearly the final part follows from the first two.

\begin{thm} Let $A$ and $Z$ satisfy hypotheses (\ref{stand}) and suppose
that $k$ is an uncountable field. Suppose that $A$ is
Cohen-Macaulay. Let $I$ be the defining ideal of the singular
locus of $Z$.

1. If $Z$ is unruffled in $A$ at the smooth points of $\Z$ then
$A[c^{-1}]$ is a flat $Z[c^{-1}]-$module for all non-zero elements
$c$ of $I.$

2. If $\m$ is a smooth point of $\Z$ such that $\m A \neq A$ and
$A_{\m}$ is a flat $Z_{\m}$-module, then $Z$ is unruffled in $A$
at $\m$.

3. Suppose that $Z$ is smooth and that $\m A \neq A$ for maximal
ideals $\m$ of $Z$. Then $A$ is a flat $Z$-module if and only if
$Z$ is unruffled in $A$.
\end{thm}

\subsection{Proof of (\ref{statement})1:} \label{proof1} Let $\m$
be a smooth point of $\Z$. We claim that, for all $i \geq 0$,
\begin{equation} \mathrm{Tor}^i_Z (Z/\m, A) \eq 0. \label{torzero}
\end{equation} By \cite[Corollary 2]{SmZh},
\begin{equation} \gk (A) \quad \geq \quad
\mathrm{GK-dim}_{Q(Z)}(Q(Z) \otimes _Z A) + \gk (Q(Z)).
\label{pjs} \end{equation} Now the unruffledness of $\m$ coupled
with (\ref{pjs}) yields
\begin{equation} \gk (A) - \gk (A/\m A) \quad \geq \quad
d.\label{gineq}
\end{equation}
Since $A$ is Cohen-Macaulay, (\ref{gineq}) implies that
\begin{equation} j(A/\m A) \quad \geq \quad d.
\label{grgrade}\end{equation}

Now Proposition (\ref{extdepth}) shows that there exist elements
$x_1, \ldots , x_d$ in $\m$ forming a regular sequence in $A$. Set
$I = \Sigma_{i=1}^d x_i Z \subseteq \m$, so that, again by
Proposition (\ref{extdepth}), \begin{equation} j(A/IA) \eq d.
\label{Igrade} \end{equation} We claim that
\begin{equation} \m \textit{ is minimal over } I. \label{dumb}
\end{equation}
For suppose (\ref{dumb}) is false, and let $\p$ be a prime of $Z$
strictly contained in $\m$ with $I \subseteq \p$, so that $\p$ has
height $r$ with $r < d.$ Then
\begin{equation} \gk (A/I A) \geq \gk (A/\p A) \geq n - r \geq n -
d + 1, \label{dumber}\end{equation} where the second inequality is
given by Lemma (\ref{primeunruff}). But (\ref{Igrade}) and
(\ref{dumber}) contradict the fact that $A$ is Cohen-Macaulay, so
(\ref{dumb}) is true. Localise in $Z$ at $\m$, so $\m^tZ_{\m}
\subseteq I_{\m}$ for some $t \geq 1$. By \cite[Corollaries 17.7
and 17.8(a)]{E} $x_1, \ldots , x_d$ constitute a $Z_{\m}$-sequence
in $Z_{\m}$ since these $d$ elements generate an ideal of the
local ring $Z_{\m}$ containing a $Z_{\m}$-sequence of length $d$,
namely the $t$th powers of a regular sequence generating $\m
Z_{\m}$. Thus the Koszul complex $K_{Z_{\m}}(x_1, \ldots , x_d)$
gives a free $Z_{\m}$-resolution of $Z_{\m}/I_{\m}.$ Since $x_1,
\ldots , x_d$ is a regular sequence in $A$, Proposition
(\ref{koszul}) shows that $K_{Z_{\m}}(x_1, \ldots , x_d)
\otimes_{Z_{\m}} A_{\m}$ has no homology except at the $d$th
place. In other words,
\[ \mathrm{Tor}^i_{Z_{\m}}(Z_{\m}/I_{\m}, A_{\m}) \eq 0 \]
for all $i > 0.$ Clearly this implies (\ref{torzero}). It follows
by Lemma (\ref{torlem}) that, if $c$ is any nonzero element of the
ideal defining the singular locus of $Z$, then $A[c^{-1}]$ is a
flat $Z[c^{-1}]$-module, and so (\ref{statement})1 is proved.

\subsection{Proof of (\ref{statement})2:}\label{proof2} Suppose that
$A$ is Cohen-Macaulay and let $\m$ be a smooth point of $\Z$ such
that $\m A \neq A$ and $A_{\m}$ is a flat $Z_{\m}$-module. Since
$Z_{\m}$ is regular $\m Z_{\m}$ is generated by a regular sequence
$x_1, \ldots , x_d$. Flatness of the $Z_{\m}$-module $A_{\m}$
ensures that $x_1, \ldots , x_d$ is a regular sequence generating
$\m A_{\m}$, as one can show easily using the Equational Criterion
for Flatness \cite[Corollary 6.5 and Exercise 6.7]{E}. Therefore
\begin{equation} j_{A_{\m}}(A_{\m}/\m A_{\m}) \eq d \label{grt}
\end{equation}
by Proposition (\ref{extdepth}). Since $\mathrm{Ext}^{*}_A (A/\m
A,A)$ is killed by $\m$ it follows from (\ref{grt}) that
$j_{A}(A/\m A) = d.$ Hence, since $A$ is Cohen-Macaulay,
\[ \gk (A/\m A) \eq \gk (A) - d. \]
Combining this with the inequality (\ref{pjs}) of (\ref{proof1})
yields
\[ \gk (A/\m A) \quad \geq \quad \mathrm{GK-dim}_{Q(Z)}(Q(Z)
\otimes_Z A). \] The reverse inequality is supplied by Lemma
(\ref{ruffineq}), so the proof of (\ref{statement})2 is complete.

\subsection{Examples.}\label{mainexamples} \textbf{1.} \textit{Theorem
(\ref{statement}).1 fails to hold whenever $A$ is any affine
commutative domain which is not Cohen-Macaulay.} For, given such
an algebra $A$, choose by Noether normalisation \cite[Theorem
13.3]{E} a polynomial subalgebra $Z$ over which $A$ is a finitely
generated module. So $Z$ is unruffled in $A$ by Example
(\ref{ruffex}).\textbf{1}. The well-known characterisation of
local commutative Cohen-Macaulay algebras by freeness over local
smooth subalgebras \cite[Corollary 18.17]{E} shows that there must
exist a maximal ideal $\m$ of $Z$ such that $A_{\m}$ is not a flat
$Z_{\m}$-module.

\textbf{2.} \textit{Theorem (\ref{statement}).2 fails in general
if $A$ is not Cohen-Macaulay (at least if $A$ is only one-sided
noetherian and is not semiprime)}. Take $A$, $Z$ and $\m$ as in
Example (\ref{ruffex}).\textbf{5}. Thus $A$ is a flat $Z$-module,
but, as we've already noted, $Z$ is not unruffled in $A$ at $\m$.
Notice that $A$ is not Cohen-Macaulay: $j(A/\m A) + \gk (A/\m A) =
1 + 0= 1 < 2 = \gk (A).$

\section{Applications}\label{apps}

\subsection{The Smith-Zhang inequality.} \label{smithy} As noted in \cite{SmZh},
the inequality (\ref{pjs}), which is their Corollary 2, is in
general strict; in fact Example (\ref{Heis}) is a case where
equality fails.\footnote{For the special case where $A$ is a
factor of an enveloping algebra, the inequality was proved in
\cite{Sm}.} However we can deduce easily from Lemma
(\ref{ruffvalue}) that the fact that this example is not
Cohen-Macaulay is the key to the failure of the equality in this
case:

\begin{cor}
Let $A$ and $Z$ be as in (\ref{stand}), and suppose that $A$ is
Cohen-Macaulay. Suppose also that $Z$ has at least one smooth
 maximal ideal for which $\m A \neq A$ and
 \begin{equation} \gk (A/\m A) \leq
\mathrm{GK-dim}_{Q(Z)}(A \otimes_Z Q(Z)). \label{dufus}
\end{equation}For example, if $A$ is finitely
related and $k$ is uncountable then this will be the case by Lemma
(\ref{constant}). Then
\begin{equation} \gk (A) \eq
\mathrm{GK-dim}_{Q(Z)}(Q(Z) \otimes _Z A) + \gk (Q(Z)).
\label{pjseq} \end{equation}
\end{cor}
\begin{proof}
This is immediate from Lemma (\ref{ruffvalue}): for this shows
that (\ref{dufus}) is an equality, both sides being equal to $n -
d.$ This proves (\ref{pjseq}).
\end{proof}

\begin{rem} The hypothesis that $A$ is Cohen-Macaulay in
Corollary (\ref{smithy}) can be relaxed a little: it is only
necessary to assume that there is a Cohen-Macaulay factor $A'$ of
$A$ with $\gk (A') = \gk (A)$, with the images of the nonzero
elements of $Z$ not zero divisors in $A'$. The adaptations needed
to the above argument are obvious.

\end{rem}

\subsection{Generalised Kostant theorem:} \label{kostant} Let $A = U(\mathfrak{g})$ be the
enveloping algebra of a finite dimensional complex Lie algebra
$\mathfrak{g}$, and let $Z$ be the centre of $A$, with $\Z =
\mathrm{maxspec}(Z)$ as usual. Example (\ref{ruffex}).4 shows that
it's not always true that $A$ is a flat $Z$-module, even when $Z$
is a polynomial algebra, as one checks in this case by direct
calculation or by appealing to Theorem (\ref{statement}).2. To
state an extra condition needed to ensure flatness, recall that
$\mathfrak{g}$ acts on the symmetric algebra $S = S(\mathfrak{g})$
via the adjoint action, and set $Y =
S(\mathfrak{g})^{\mathfrak{g}}$. Thus $S$ is the associated graded
algebra of the filtered $\C$-algebra $A$, and the canonical map
from $A$ to $S$ is an isomorphism of $\mathfrak{g}$-modules which
carries $Z$ to $Y$ and has as inverse the \emph{symmetrisation
map}, \cite[Proposition 2.4.10]{Dix}.

\begin{thm} Retain the above notation, and assume that $Z$ is
affine and that $\m A \neq A$ for all maximal ideals $\m$ of $Z$. Let $\mathfrak{y}^{+}$
be the augmentation ideal of $Y$, that is $\mathfrak{y}^{+} =
\mathfrak{g}S \cap Y.$

1.  Consider the statements: \begin{enumerate} \item
$\mathfrak{y}^{+}$ is a smooth point of $Y$ and is unruffled in
$S$. \item $S_{\mathfrak{y}^{+}}$ is a flat
$Y_{\mathfrak{y}^{+}}$-module. \item $Z$ is unruffled in $A$ at
$\m$ for all smooth points $\m$ of $\Z$. \item $A_{\m}$ is a flat
$Z_{\m}$-module for all smooth points $\m$ of $\Z$.
\end{enumerate} Then $$ (1)\Leftrightarrow(2)\Rightarrow (3) \Leftrightarrow (4).$$

2. Suppose in addition that $Z$ is smooth (which is equivalent to
assuming that $Z$ is a polynomial algebra). If $\mathfrak{y}^{+}$
is unruffled in $S$ then $A$ is a flat $Z$-module.

\end{thm}

\begin{proof} 1. By \cite[Theorem 10.4.5]{Dix}, $Y$ and $Z$ are isomorphic
 (although not in general via the symmetrisation map); in particular, $Y$
 is affine since we are assuming that $Z$ is. Since $S$ is a polynomial algebra and so in
particular smooth, $\mathfrak{y}^{+}$ is a smooth point of $Y$ if
$S_{\mathfrak{y}^{+}}$ is a flat $Y_{\mathfrak{y}^{+}}$-module,
since a finite $S_{\mathfrak{y}^{+}}$-projective resolution of the
trivial $S$-module yields a finite flat resolution of the unique
simple $Y_{\mathfrak{y}^{+}}$-module. Thus the equivalence of (1)
and (2) follows from the commutative case of the Main Theorem
(\ref{statement}); see (\ref{ruffex}).\textbf{3}.

Suppose now that (1) and hence (2) hold. Since the associated
graded algebra $S$ of $A$ is smooth, $A$ is Cohen-Macaulay by
\cite[Theorem II.2.1 ]{Bj}. We claim that \begin{equation}
\mathrm{GK-dim}_{Q(Z)}(Q(Z) \otimes_Z A) \eq
\mathrm{GK-dim}_{Q(Y)}(Q(Y) \otimes_Y S). \label{gradeeq}
\end{equation}By Corollary (\ref{smithy}),
\begin{equation}
\mathrm{GK-dim}_{Q(Z)}(Q(Z) \otimes_Z A) \eq \gk (A) - \gk (Z),
\label{smart}
\end{equation} and similarly (although in this case \cite[Theorem
13.5]{E} suffices),\begin{equation} \mathrm{GK-dim}_{Q(Y)}(Q(Y)
\otimes_Y S) \eq \gk (S) - \gk (Y). \label{popart} \end{equation}
But $S$ and $Y$ are respectively the associated graded algebras of
$A$ and $Z$, so the right hand sides of (\ref{smart}) and
(\ref{popart}) are equal by \cite[Proposition 6.6]{KL}, proving
our claim.

Next, as $S_{\mathfrak{y}^{+}}$ is a flat
$Y_{\mathfrak{y}^{+}}$-module, Lemma (\ref{ruffineq}) implies that
\begin{equation}\mathrm{GK-dim}_{Q(Y)}(Q(Y) \otimes_Y S) \geq \gk(S/
\mathfrak{y}^{+}S).\label{sink} \end{equation} Now let $\m^{+} =
\mathfrak{g}A \cap Z,$ the augmentation ideal of $Z$. Let $\m$ be
a smooth point of $\Z$. Thus, writing $\mathrm{gr}(-)$ for
associated graded modules,
\begin{equation}\mathrm{gr}(\m A) \quad \supseteq \quad \mathrm{gr}(\m)S \eq
\mathrm{gr}(\m^{+})S \eq \mathfrak{y}^{+}S.
\label{graft}\end{equation} Hence, \begin{equation}
\gk(S/\mathfrak{y}^{+}S) \geq \gk (A/\m A). \label{grafter}
\end{equation} By (\ref{gradeeq}), (\ref{sink}) and
(\ref{grafter}), \begin{equation} \mathrm{GK-dim}_{Q(Z)}(A
\otimes_Z Q(Z)) \geq \gk (A/\m A). \label{push} \end{equation} But
$\m A \neq A,$ so Lemma (\ref{ruffineq}) applies and $Z$ is
unruffled in $A$ at $\m$. That is, (2) $\Rightarrow$ (3). The
equivalence of (3) and (4) follows from the Main Theorem
(\ref{statement}).

 2. Suppose now that $Z$ is smooth. Thus so
also is $Y$ by \cite[Theorem 10.4.5]{Dix}. So the result follows
from (1) $\Rightarrow$ (4) of the first part. \end{proof}

\begin{cor}($\mathrm{Kostant},$ \cite[Theorem 8.2.4]{Dix})
Suppose that $\mathfrak{g}$ is a finite dimensional complex
semisimple Lie algebra. Let $A = U(\g)$ and let $Z$ be the centre
of $A$. Then $A$ is a free $Z$-module.
\end{cor}
\begin{proof} Retain the notation of the theorem. We check that the
hypotheses of the second part of the theorem are satisfied. When
$\mathfrak{g}$ is semisimple $Z$ is a polynomial algebra on
$\mathrm{rank}(\mathfrak{g})$ indeterminates, \cite[Theorem
7.3.8(ii)]{Dix}. Local finiteness of the adjoint action of
$\mathfrak{g}$ on $A$ combined with the semisimplicity of finite
dimensional $\mathfrak{g}$-modules imply that $A$ is a direct sum
of finitely generated $Z$-modules, and so $\m A \neq A$ for each
maximal ideal $\m$ of $Z$, by Nakayama's Lemma. The subvariety in
$\mathfrak{g}^* = \mathfrak{g}$ defined by $\mathfrak{y}^+
S(\mathfrak{g})$ is the cone of nilpotent elements, \cite[Theorem
8.1.3(i) ]{Dix}, which has dimension
$\mathrm{dim}_{\C}(\mathfrak{g}) - \mathrm{rank}(\mathfrak{g})$ by
\cite[Theorem 8.1.3(ii)]{Dix}. So $\mathfrak{y}^+ S(\mathfrak{g})$
is unruffled in $S(\mathfrak{g})$. Thus $A$ is a flat $Z$-module
by the second part of the theorem.

Since $A$ is a direct sum of finitely generated (and so
projective) $Z-$modules, and projective modules over the
polynomial algebra $Z$ are (stably) free, flatness implies
freeness in this case.
\end{proof}

\subsection{Questions of Borho and Joseph} \label{questions} Let $\g$
be a finite dimensional complex semisimple Lie algebra. In a
series of papers \cite{BJ,Bo2,Bo3,Bo4} Borho and Joseph have
studied the primitive spectrum $\chi$ of $U(\g)$ by partitioning
$\chi$ into sheets. By definition, a \emph{sheet} in $\chi$ is an
irreducible subset $\mathcal{Y}$ of $\chi$ which is maximal such
that $\gk (U(\g)/P)$ is constant for $P \in \mathcal{Y}$ and the
Goldie dimension of $U(\g)/P$ is bounded for $P \in \mathcal{Y}$.
In \cite[Corollary 5.6]{BJ} it is shown that every sheet in $\chi$
has the form $\overline{\chi}(J,\mathfrak{z})$, where the latter
is defined as follows.

Let $\mathfrak{h}$ be a Cartan subalgebra of $\g$ and let
$\mathfrak{p}$ be a parabolic subalgebra of $\g$ with
$\mathfrak{h} \subseteq \p$ and with Levi decomposition
$\mathfrak{p} = \mathfrak{m} \oplus \mathfrak{l}$, and let
$\mathfrak{z}$ be the centre of $\mathfrak{l}$. Let $J$ be a
primitive completely rigid ideal of $U(\mathfrak{l})$; this means
that $J$ is not almost induced\footnote{A primitive ideal of
$U(\mathfrak{l}))$ is almost induced if it is a minimal prime over
an ideal of the form $I_{\mathfrak{p'}}(J', \mu)$ for a parabolic
subalgebra $\p'$ of $\mathfrak{l}$.} from any proper Levi
subalgebra of $\mathfrak{l}$. (See \cite[5.6]{BJ} for details; for
example, a primitive ideal of finite codimension is completely
rigid, but not conversely in general.) For $\lambda \in
\mathfrak{z}^*$, define $I_{\mathfrak{p}}(J, \lambda)$ to be the
annihilator in $U(\g)$ of $U(\g) \otimes_{U(\mathfrak{p})}
((U(\mathfrak{l})/J) \otimes \C_{\lambda})$, where $\C_{\lambda}$
denotes the one-dimensional $U(\mathfrak{p})-$module with weight
$\lambda$, where we identify $\mathfrak{z}$ with $\p/[\p,\p]$ in
order to view $U(\mathfrak{z})-$modules as $U(\p)-$modules. Then
\begin{equation} \overline{\chi}(J,\mathfrak{z}) \eq \{ I \in \chi : I
\textit{ is minimal over } I_{\mathfrak{p}}(J, \lambda), \quad
\lambda \in \mathfrak{z}^* \}. \label{sheet}\end{equation} Another
way of describing $\overline{\chi}(J,\mathfrak{z})$ is as the set
of minimal primitive ideals of the prime factor $A = U(\g)/P$ of
$U(\g)$, where
$$ P \eq \cap_{\lambda \in \mathfrak{z}^*}I_{\mathfrak{p}}(J, \lambda).$$
Fix a weight $\nu$ such that $J$ is the annihilator of the
irreducible highest weight $U(\mathfrak{l})-$module $L'(\nu)$.
Here, we can take $\nu \in \mathfrak{z}^{\perp}$, the Killing
orthogonal to $\mathfrak{z}$ in $\mathfrak{h}$, so that
$\mathfrak{z}^{\perp}$ is a Cartan subalgebra of
$[\mathfrak{l},\mathfrak{l}]$. Then $P$ is the annihilator in
$U(\mathfrak{g})$ of $U(\g) \otimes_{U(\p)} (L'(\nu)\otimes
U(\mathfrak{z}))$.

Thus, to study the sheets in $\chi$ amounts to studying the
collection of minimal primitive ideals of the factors of $U(\g)$
of the form $A$. In particular, with the notation we've introduced
above, the sheet $\overline{\chi}(J,\mathfrak{z})$ consists
precisely of (the inverse images in $U(\g)$ of) the prime ideals
of $A$ which are minimal over an ideal generated by a maximal
ideal of $Z$, the centre of $A$. As is implied by Proposition
(\ref{generic}).4, for a dense set of those maximal ideals $\m$ of
$Z$, $\m A$ is in fact prime and hence primitive. However
typically there are exceptional $\m$ for which this is not the
case, and in an attempt to remedy this one passes to the larger
algebra
$$\tilde{A} \quad := \quad A \otimes_Z \tilde{Z}, $$ where
$\tilde{Z}$ is the integral closure of $Z$ in its quotient field.
It is still not always true that $I = (I \cap \tilde{Z})\tilde{A}$
for every minimal primitive ideal $I$ of $\tilde{A}$
\cite[4.6]{Bo2}, but Borho proves in \cite[\S 9, Theorem]{Bo3}
that, at least when J is the augmentation ideal of
$U(\mathfrak{l})$, every minimal primitive $I$ of $\tilde{A}$
satisfies
$$ I \eq \sqrt{(I \cap \tilde{Z})\tilde{A}}. $$

The analysis of $\tilde{A}$ and $\overline{\chi}(J,\mathfrak{z})$
would be greatly facilitated if a positive answer to the following
question from \cite{J} , were known. (See also the closely related
question in \cite[5.3,
Remark (b)]{BJ}.)\\

\noindent \textbf{Question 1:} Is $\tilde{A}$ a free $\tilde{Z}-$module?\\

As we've noted in Example (\ref{ruffex})\textbf{2}, $A$ is
unruffled, and the same argument from \cite[5.8]{BJ} shows that
\begin{equation} \tilde{A} \textit{ is unruffled over }\tilde{Z}.
\label{tilunruff}
\end{equation}
Thus it's clear from Theorem (\ref{statement}) that
Question 1  is closely connected to \\

\noindent \textbf{Question 2:} With the above notation, is
$\tilde{A}$
Cohen-Macaulay?\\

We don't know the answer to this question. We shall show here
however that, at least in an important special case, a positive
answer to Question 2 implies a positive answer to Question 1.
Retain all the notation already introduced in this subsection. Let
$\hat{W}$ be the normaliser of $\mathfrak{z}$ in the Weyl group
$W$ of $\g$, and let $\hat{W}_{\nu} = \{w \in \hat{W}: w\nu = \nu
\}.$ By \cite[Proposition 6.1b]{Bo5} or \cite[proof of Proposition
8.6(b)]{BJ},
$$ \tilde{Z} \eq S(\mathfrak{z})^{(\hat{W}_{\nu},*)}, $$
where $*$ denotes the shifted action, $w*\lambda = w(\lambda +
\rho') - \rho',$ and where \linebreak $\rho' =
-\frac{1}{2}(\textit{sum of roots in } \g/\p)\mid_{\mathfrak{z}}$.
Now assume that \begin{equation} \hat{W}_{\nu} \textit{ is
generated by reflections,} \label{reflect}
\end{equation}so that, by the
Shepherd-Todd-Chevalley theorem,
\begin{equation} \tilde{Z} \textit{ is a polynomial algebra. }
\label{centre}\end{equation}

\begin{thm} Retain the notation introduced in this subsection.
Assume (\ref{reflect}). If $\tilde{A}$ is Cohen-Macaulay, then
$\tilde{A}$ is a free $ \tilde{Z}-$module.
\end{thm}

\begin{proof} Assume (\ref{reflect}) and that $\tilde{A}$ is
Cohen-Macaulay. In view of (\ref{centre}) and (\ref{tilunruff}),
the hypotheses of Theorem (\ref{statement}).3 are satisfied, so we
can conclude that $\tilde{A}$ is a flat $\tilde{Z}-$module. Thanks
to the local finiteness and complete reducibility of the adjoint
action of $\g$ on $U(\g)$, $\tilde{A}$ is a direct sum of finitely
generated $\tilde{Z}-$modules, so that $\tilde{A}$ is a free
$\tilde{Z}-$module as claimed.
\end{proof}

\begin{rems} 1. When $\g = \mathfrak{sl}(n)$, $\hat{W}$ is always
generated by reflections. Moreover the completely rigid primitive
ideal $J$ of $U(\mathfrak{l})$ will always in the
$\mathfrak{sl}(n)$ case be co-artinian, \cite[6.10]{BJ}. Thus if
we are concerned only with sheets of completely prime primitive
ideals in $U(\mathfrak{sl}(n))$, then $J$ will always be the
augmentation ideal of $U(\mathfrak{l})$, so $\nu = 0$, and
(\ref{reflect}) is satisfied.

2. There are some tentative indications that ``many'' prime
factors of enveloping algebras $U(\g)$ of semisimple Lie algebras
may be Auslander-Gorenstein\footnote{The definition is recalled in
(\ref{factors}).} and/or Cohen-Macaulay. For example, if $P$ is a
maximal ideal of $U(\g)$ then $U(\g)/P$ is Auslander-Gorenstein by
\cite{Y}. On the other hand, if $P$ is a minimal primitive ideal
then the same conclusion holds by \cite{Lev2}. This latter result
can be generalised: if $P$ is any primitive ideal of $U(\g)$ for
which (a) gr$(P)$ is prime and (b) the closure
$\overline{\mathcal{O}}$ of the associated (nilpotent) orbit
$\mathcal{O}$ of $P$ is Gorenstein, then standard filtered-graded
arguments yield that $U(\g)/P$ is Auslander-Gorenstein and
Cohen-Macaulay. Sufficient conditions for (a) to hold can be read
off from \cite{BB}; and (b) always holds for the normalisation of
$\overline{\mathcal{O}}$ by \cite{H} or \cite{Pan}, and hence
always holds for $\overline{\mathcal{O}}$ itself in type $A$,
\cite{HP}. One can then hope to lift such a property from
$U(\g)/P$ to the closure of the sheet containing $P$. But since we
have only partial results in this direction we shall not pursue
this here.
\end{rems}

\subsection{Factors of unruffled algebras} \label{factors} Suppose
that $Z$ and $A$ satisfy (\ref{stand}),
with $Z$ unruffled in $A$, and $P$ is (say) a prime ideal of $A$.
Is $A/P$ unruffled over $Z/P \cap Z$? The example below shows that
the answer is no even when $A$ is commutative. However it may be
that positive results can be obtained in special circumstances;
for instance, this might be one route to a more elementary proof
of \cite[Corollary 5.8]{BJ}, that prime factors of the enveloping
algebra of a complex semisimple Lie algebra $\g$ are unruffled
over their centres, since it is relatively easy to see that
$U(\g)$ itself is unruffled over its centre, Theorem
(\ref{kostant}).2. In this subsection we show that, at least under
some additional hypotheses, the unruffled property is stable under
factoring by an ideal of $A$ generated by a prime ideal of $Z$. On
the way we derive some useful subsidiary
results.  \\

\noindent \textbf{Example:} Let $A = \C [x,y,z]$ and let $Z$ be
the subalgebra $\C [x,yz]$ of $A$. It's trivial to check that $Z$
is unruffled in $A$; equivalently (by Theorem (\ref{statement}).3,
$A$ is a flat $Z$-module. But if we factor by the prime ideal
$(x-z)A$ we
get the ruffled example (\ref{ruffex}).\textbf{3}.   \\

\begin{thm}  Assume that $A$ and $Z$ satisfy (\ref{stand}) and (\ref{uncount})
of (\ref{generic}). Suppose that $Z$ is smooth, and unruffled in
$A$. Let $\p$  be a prime ideal of $Z$ of height $\ell$. Then:
\begin{enumerate} \item $Z/\p$ is unruffled in
$A/\p A.$\\
Suppose in addition that $A$ is Cohen-Macaulay. Let $F$ denote the
field of fractions of $Z/\p .$ Then: \item $\mathrm{GK-dim}_F
(A/\p A \otimes_{Z/\p} F) \eq n - d;$
\item $\gk (A/\p A) \eq n - \ell.$
\end{enumerate}
\end{thm}

\begin{proof} The Main Theorem (\ref{statement}).3 implies
that $A$ is a flat $Z-$module. Hence, by the Equational Criterion
for Flatness,
 \cite[Corollary 6.5 and Exercise 6.7]{E}, \begin{equation} \textit{the elements
 of }
 Z \setminus \p \textit{ are not zero divisors in } A/\p A. \label{zero}\end{equation}
 In particular, $Z/\p \subseteq A/\p A$,
 and this pair of algebras satisfies the hypotheses (\ref{stand}) and (\ref{uncount}) of (\ref{generic}).

 1. Since $Z$ is unruffled in $A$, the GK-dimension of the factors
 $(A/\p A)/(\m A/\p A)$ is constant as $\m$ ranges through the
 maximal ideals of $Z$ which contain $\p.$ Hence unruffledness of
 $Z/\p$ in $A/\p A$ follows from Lemma (\ref{constant}).

 2. This is immediate from 1. and Lemma (\ref{ruffvalue}).

 3. The desired result is true when $\p = 0$ and also, by 2., when
 $\p$ is a maximal ideal of $Z$. Since $Z$ is an affine domain
 there is a chain $0 = \p_0 \subset \p_1 \subset \ldots
 \subset \p_{\ell} = \p \subset \dots\subset \p_d$ of
 prime ideals of $Z$. By (\ref{zero}) and \cite[Proposition 3.15]{KL}, the
 GK-dimension of the factors $A/\p_i A$ goes down by at least one
 as we pass up each step of the chain. Since the difference
 between the GK-dimensions of $A$ and of $A/\p_d A$ is exactly
 $d$, the GK-dimension must go down by exactly one at each step,
 as required.
\end{proof}

Parts (2) and (3) of the theorem fail without the Cohen-Macaulay
hypothesis, as is shown by Example (\ref{primeunruff}).

We can improve (3) of the theorem by showing that $A/\p A$ is
GK-homogeneous, but only under the extra - presumably superfluous
- hypothesis that $A$ is Auslander-Gorenstein. Recall that a
Noetherian ring $R$ is {\em Auslander-Gorenstein} if the
$R$-module $R$ has finite (equal) right and left injective
dimensions, and $R$ satisfies the Auslander conditions; namely,
for every non-zero left or right $R$-module $M$ and every
non-negative integer $i$, every non-zero submodule $N$ of
$\textrm{Ext}_R^i(M,R)$ satisfies $j_R(N) \geq i$. Details and
further references can be found in \cite{Lev}, for example.

\begin{lem}
Let $R$ be a Noetherian Auslander-Gorenstein $k-$algebra, let $z$
be a central element of $R$, and let $V$ be a non-zero finitely
generated $R-$module on which $z$ acts torsion freely. Then $j_R
(V) = j_{R[z^{-1}]}(V \otimes_R R[z^{-1}]).$
\end{lem}

\begin{proof}
It's clear from the behaviour of $\mathrm{Ext-}$groups under
central localisation that $j_R (V) \leq j_{R[z^{-1}]}(V \otimes_R
R[z^{-1}]).$ To prove the reverse inequality, set $\overline{V} =
V/Vz$, so that there is an exact sequence
$$ 0 \longrightarrow V \overset{z \times}\longrightarrow V \longrightarrow
\overline{V}\longrightarrow 0. $$ The part of the long exact
sequence of $\textrm{Ext}$-groups around $j := j_R(V)$ is thus
$$  \textrm{Ext}_R^{j}(\overline{V},R) \To \textrm{Ext}_R^{j}(V,R)
 \overset{z \times } \To \textrm{Ext}_R^{j}(V,R)\To
\textrm{Ext}_R^{j+1}(\overline{V},R).$$ Here,
$\textrm{Ext}_R^{j}(\overline{V},R) = 0$ by \cite[Theorem
4.3]{Lev}, since $R$ is Auslander-Gorenstein, showing that
$\textrm{Ext}_R^{j}(V,R)$ has no $\{z^i\}-$torsion. Thus $j_R (V)
\geq j_{R[z^{-1}]}(V \otimes_R R[z^{-1}]),$ as required.
\end{proof}

\begin{prop} Assume that $A$ and $Z$ satisfy (\ref{stand}) and (\ref{uncount})
of (\ref{generic}). Suppose that $A$ is Auslander-Gorenstein and
Cohen-Macaulay, and that $Z$ is smooth, and unruffled in $A$. Let
$\p$ be a prime ideal of $Z$ of height $\ell.$ Then there exists
an element $z \in Z\setminus \p$ such that $(A/\p A)[z^{-1}]$ is
Auslander-Gorenstein and Cohen-Macaulay of dimension $n - \ell.$.
\end{prop}
\begin{proof}
Since $Z$ is smooth, there exists an element $z$ in $Z \setminus
\p$ such that $\p[z^{-1}]$ is generated by a regular sequence
$\{x_1, \ldots , x_{\ell}\}$ in $Z$. As in the proof of the lemma,
$\{x_1, \ldots , x_{\ell} \}$ is a regular sequence in $A$. Thus
$(A/\p A)[z^{-1}]$ is Auslander-Gorenstein by \cite[3.4, Remark
(3)]{Lev}.\\
To prove the Cohen-Macaulay property, let $L$ be a finitely
generated $(A/\p A)[z^{-1}]-$module. By \cite[Corollary 11.68]{R}
\begin{equation} j_{(A/\p A)[z^{-1}]}(L) \eq j_{A[z^{-1}]}(L) -
\ell. \label{rot} \end{equation} Fix a finitely generated
$A-$submodule $L_0$ of $L$ such that $L = L_0 \otimes A[z^{-1}].$
Then
\begin{equation} j_{A[z^{-1}]}(L) \eq j_{A}(L_0)
\label{ore}\end{equation} by the Lemma above and, since $A$ is
Cohen-Macaulay,
\begin{equation} j_A(L_0) \eq n - \gk (L_0). \label{macked}
\end{equation}
Since $\gk (L_0) = \gk (L)$ by \cite[Proposition 4.2]{KL}, from
(\ref{rot}), (\ref{ore}) and (\ref{macked}) it follows that
\begin{equation} j_{(A/\p A)[z^{-1}]}(L) \eq n - \ell - \gk (L),
\end{equation}
and this combined with Theorem (\ref{factors}) proves the result.
\end{proof}

\begin{cor}Assume that $A$ and $Z$ satisfy (\ref{stand}) and (\ref{uncount})
of (\ref{generic}). Suppose that $A$ is Auslander-Gorenstein and
Cohen-Macaulay, and that $Z$ is smooth, and unruffled in $A$. Let
$\p$ be a prime ideal of $Z$ of height $\ell.$ Then $A/\p A$ is
GK-homogeneous of dimension $n - \ell$, and has an Artinian
quotient ring.
\end{cor}
\begin{proof}
Let $z \in Z \setminus \p$ be the element afforded by the
proposition. By (\ref{zero}), $A/\p A$ embeds in $(A/\p
A)[z^{-1}]$, and it is easy to check by a small adjustment to the
proof of \cite[Proposition 4.2]{KL} that it is enough to prove
that the desired conclusions hold for $A/\p A[z^{-1}].$ Now the
case of grade zero of the Cohen-Macaulay property implies
GK-homogeneity of $(A/\p A)[z^{-1}]$. That this implies the
existence of an Artinian quotient ring for $A/\p A$ now follows
from \cite[Theorem 5.3]{Lev}.
\end{proof}

\section{Questions}\label{qns} Some of the questions listed here
have already been mentioned earlier; we record them again for the
reader's convenience.
\subsection{GK-dimension} \label{GKdim} Is there a generalisation of the
Main Theorem to settings where GK-dimension is not defined? In
particular, is there a good way to define the Cohen-Macaulay
property in the absence of GK-dimension?

Example (\ref{Heis}), the Heisenberg group algebra $kH$, is
Cohen-Macaulay neither with the definition (\ref{cohmac}) used in
this paper, nor with the definition using the Krull dimension.
Moreover, if one defines {\em Krull unruffled} extensions in the
obvious way, using the Krull dimension rather than the
GK-dimension, then $kH$ is {\em not} Krull unruffled over its
centre $Z$. Nevertheless, $kH$ is free over $Z$, which of course
is smooth and affine. Is there a version of the Main Theorem
incorporating this example? For example, perhaps the correct
setting is that of algebras $A$ for which there is an integer
$\mu$ such that $\delta(-) := \mu - j_A(-)$ defines an exact
finitely partitive dimension function? By \cite[Definition
4.5]{Lev} this would include Auslander-Gorenstein algebras such as
$kH$.

\subsection{Density of unruffled points} \label{dense} (Borho-Joseph, \cite{BJ}; see Lemma (\ref{constant}).) Suppose
that $A$ and $Z$ satisfy (\ref{stand}) and (\ref{uncount}) of
(\ref{generic}). Does the set of unruffled maximal ideals of $Z$
contain a non-empty Zariski open subset of $\Z$?

\subsection{Existence of unruffled extensions} \label{exist} Find
a more elementary proof of the fact that prime factors of the
enveloping algebra of a complex semisimple Lie algebra are
unruffled over their centres. Find some other large classes of
unruffled extensions. For example, what about quantised enveloping
algebras $U_q(\g)$ where the quantising parameter $q$ is not a
root of unity? Is a prime noetherian affine PI algebra with an
affine centre $Z$ always unruffled over $Z$?

\subsection{Non-central subalgebras} \label{noncentral} Let $B
\subseteq A$ be any pair of affine noetherian algebras of finite
GK-dimension. One can clearly extend the definition of an
unruffled extension to this setting: several variants come to
mind, but one might try requiring constancy of $\gk (M \otimes_B
A)$ as $M$ ranges over all simple right $B-$modules of fixed
GK-dimension. Does this lead to an interesting theory? Is there a
version of the Main Theorem?

\subsection{Unruffled factors} \label{factor2} Is Corollary
(\ref{factors}) true without the hypothesis that $A$ is
Auslander-Gorenstein? Are all the results of this paragraph valid
without assuming $Z$ smooth, provided $\p$ is a smooth prime?

\subsection{Factors of semisimple enveloping algebras}
\label{semisimple} First, we repeat (a generalisation of) Joseph's
question from \cite{J}, already stated in (\ref{questions}). Let
$P$ be a prime ideal of the enveloping algebra of a complex
semisimple Lie algebra $\g$, and suppose $P$ is induced from a
completely rigid primitive ideal $J$ of the enveloping algebra of
a parabolic subalgebra $\p$, $P = I_{\p}(J,\lambda).$ Let
$\tilde{Z}$ be the normalisation of the centre $Z$ of $A =
U(\g)/P.$ Suppose that $\tilde{Z}$ is smooth. Is $\tilde{A} := A
\otimes_Z \tilde{Z}$ a free $\tilde{Z}-$module?

In view of Theorem (\ref{questions}) a positive answer to the
above question would follow from a positive answer to the
following. With the above notation and hypotheses, is $\tilde{A}$
Cohen-Macaulay? One can also ask, of course, whether $\tilde{A}$
is Auslander-Gorenstein.

More generally, the partial results for primitive ideals discussed
in (\ref{questions}) suggest the following rather wild
speculation: if $\tilde{A}$ is the normalisation of an arbitrary
prime factor of $U(\g)$, $\g$ semisimple, is the
Auslander-Gorenstein and/or the Cohen-Macaulay property for
$\tilde{A}$ controlled by the corresponding property for
$\tilde{Z}$? In particular, which primitive factors of $U(\g)$ are
Auslander-Gorenstein or Cohen-Macaulay?

% ----------------------------------------------------------------
\bibliographystyle{amsplain}

\end{document}